\documentclass[oneside,leqno,11pt]{article}
\usepackage{amssymb,amsmath,latexsym}

\def\gg{\mathfrak{g}}

\def\gk{\mathfrak{k}}
\def\gl{\mathfrak{l}}
\def\gm{\mathfrak{m}}

\def\gp{\mathfrak{p}}

\def\gs{\mathfrak{s}}

\def\gv{\mathfrak{v}}

\def\gz{\mathfrak{z}}

\def\C{\mathbb{C}}

\def\H{\mathbb{H}}

\def\R{\mathbb{R}}

\newtheorem{theorem}[equation]{Theorem}

\newtheorem{lemma}[equation]{Lemma}

\newtheorem{proposition}[equation]{Proposition}

\title{Complex forms of Quaternionic Symmetric Spaces}

\author{Joseph A. Wolf
\thanks{
Research partially supported by NSF Grant
DMS 99-88643.\endgraf
{\em 2000 AMS Subject Classification.} Primary 53C26;
secondary 53C29, 53C30.\endgraf
{\em Key Words}: quaternionic manifold, complex submanifold, complex form.}}

\date{19 September 2003}

\begin{document}

\maketitle

\abstract{We give a complete classification of the  complex forms of 
quaternionic symmetric spaces.}

\section{Introduction}
\setcounter{equation}{0}

Some years ago, H. A. Jaffee found the real forms of hermitian symmetric 
spaces (\cite{J1}, \cite{J2}; or see \cite{HO}).  That classification turns 
out to be related to the classification of causal symmetric spaces. 
This was first observed by I. Satake (\cite[Remark 2 on page 30]{S}
and \cite[Remark on page 87]{S}).  Somewhat later it was independently 
observed by J. Hilgert, G. \' Olafsson and B. \O rsted; see \cite{HO},
especially Chapter 3 and the Notes at the end of that Chapter.  I learned
about that from Bent \O rsted.  He and Gestur \' Olafsson had informally
discussed complex forms of quaternionic symmetric spaces and found examples
for the classical groups, for $G_2$\,, and perhaps for $F_4$\,.   \O rsted told
me about the classical ones, and we rediscovered examples for $G_2$ and $F_4$\,.
I thank Bent \O rsted for agreeing to my incorporating those examples into 
this note.  Later I used the computer program LiE \cite{L} to find examples for 
$E_6$\,, $E_7$ and $E_8$\,.
\medskip

In this note I write down a complete classification for complex
forms $L/V$ of quaternionic symmetric spaces $G/K$.  The definitions and
some preliminary results are in Sections \ref{qss} and \ref{defs}, the
main results are stated in Section \ref{classification}, and the proofs
are in Sections \ref{equal-rank-classical}, \ref{prog},  
\ref{equal-rank-exceptional} and \ref{unequal-rank}.  The case where $G$ is
a classical group and rank$(L)$ = rank$(G)$ is handled, essentially by 
matrix considerations, in Section \ref{equal-rank-classical}.  That, of course,
does not work comfortably for the exceptional groups, which must be approached
by means of their root structure.  The tool for this is a script for use of
the computer program LiE; it is described in Section \ref{prog} along with
some examples of its application.  Those examples have the interesting
property that the complexifications $L_\C$ and $K_\C$ are conjugate 
in $G_\C$\,.  They cover the delicate cases
for $G$ exceptional and  rank$(L)$ = rank$(G)$, and the remaining exceptional
equal rank cases are settled in Section \ref{equal-rank-exceptional}.  Finally
the few cases of rank$(L) <$ rank$(G)$ are worked out in Section
\ref{unequal-rank}.
\medskip

A possible extension of the theory is mentioned in Section \ref{quat-form}.
\medskip

After this paper was written I learned that quite a lot was published on
totally complex submanifolds of quaternionic symmetric spaces from the
viewpoint of differential geometry.  See, for example, \cite{ADM}, \cite{AM1},
\cite{AM2}, \cite{F}, \cite{JKS}, \cite{L1}, \cite{L2}, \cite{Ma}, \cite{Mo},
\cite{Ts} and \cite{X}, but especially the first three.  I also learned 
that M. Takeuchi \cite{Ta} had studied the maximal totally complex 
submanifolds of quaternionic symmetric spaces, reducing their classification
to that of certain Satake diagrams and writing out the classification in the
classical group cases.  {\it A priori} that is not quite the same
as the classification of complex forms of quaternionic symmetric spaces,
but it is very close.  On the other hand it seems to me that the method given
here is more efficient and more direct, and more explicit in the exceptional
group cases.  I thank Dmitry Alexeevsky for calling the above--cited papers 
to my attention.
\medskip

\section{Quaternionic Symmetric Spaces} \label{qss}
\setcounter{equation}{0}

We recall the structure of quaternionic symmetric spaces \cite{W}
A {\bf quaternionic structure} on a connected riemannian manifold $M$ is a 
parallel field $A$ of quaternion algebras $A_x$ on the real tangent spaces 
$T_x(M)$\,, such
that every unimodular element of every $A_x$ is an orthogonal linear
transformation.  Thus $A$ gives every tangent space the structure of 
quaternionic
vector space such that the riemannian metric at $x$ is hermitian relative to the elements of $A_x$ of square $-I$.  If $n = \dim M$ then a quaternionic structure 
is the same as a reduction of the structure group of the tangent bundle from $O(n)$
to $Sp(n/4)\cdot Sp(1)$  Let $K_x$ denote the holonomy group of $M$ at $x$
(we will see in a minute that this is appropriate notation for symmetric spaces
with no euclidean factor).  Suppose that $M$ is simply connected, so that
the $K_x$ are connected.  Let $A = \{A_x\}$ be a quaternionic structure on $M$.
Then $A_x$ is stable under the action of $K_x$\,, so $K_x \cap A_x$ is a
closed normal subgroup of $K_x$\,.  Now $K_x = K_x^{lin}\cdot K_x^{sca}$ where
$K_x^{lin}$ is the quaternion--linear part, centralizer of $A_x$ in $K_x$\,, and$K_x^{sca} = K_x \cap A_x$ is the scalar part.  We say that
$K_x$ has {\bf real scalar part} if $K_x^{sca}$ consists of real scalars,
i.e. $K_x^{sca}$ is $\{1\}$ or $\{\pm 1\}$.  We say that $K_x$ has {\bf complex
scalar part} if $K_x^{sca}$ is contained in a complex subfield of $A_x$ but
not in the real subfield, and we say that $K_x$ has {\bf quaternion scalar part} 
if $K_x^{sca}$ is not contained in a complex subfield of $A_x$\,.

A riemannian $4$--manifold $M$ with holonomy $U(2)$ has a dual role: it has a
quaternionic structure $A_1$ generated by the $SU(2)$--factor in the holonomy;
that has quaternionic scalar part, the same $SU(2)$,
$M$ it has a second quaternionic structure 
$A_2$ where $A_{2,x}$ is the centralizer of $A_{1,x}$ in the
algebra of $\mathbb R$--linear transformations of $T_x(M)$;  it has
complex scalar part, generated by the circle center of the holonomy $U(2)$.
Thus we have an interesting dual picture: the holonomy of $M$ has quaternionic 
scalar part for $A_1$ and has complex scalar part for $A_2$\,.
 
\begin{proposition} \label{quaternion-criterion}
The connected simply connected riemannian symmetric spaces with quaternionic
structure are the following.
 
{\rm (i)} The euclidean spaces of dimension divisible by $4$.  
Here the holonomy has real scalar part.
 
{\rm (ii)} Products $M = M_1 \times \dots \times M_\ell$ where each $M_i$ is
(a) the complex projective
or hyperbolic plane with the quaternionic structure of complex scalar part, or
(b) a product $M_i' \times M_i^"$ where each factor is a complex
projective line a complex hyperbolic line.  Here $M = G/K$, $K$ is the
holonomy, and the holonomy has complex scalar part.
 
{\rm (iii)} Irreducible connected simply connected riemannian symmetric spaces
$M = G/K$ where $K$ has an $Sp(1)$ factor that generates  quaternion algebras
on the tangent spaces of $M$.  Here $K$ is the holonomy, and the holonomy has
quaternion scalar part.
\end{proposition}
 
There is a structure theory for the spaces of Proposition
\ref{quaternion-criterion}(iii).  There are two, a compact one and its noncompact
dual, for each complex simple Lie algebra, and they are constructed from the
highest root \cite{W}.  These spaces are listed in the Table
\ref{quaternion-symm} below.  Here we use the notation that $G_2$\,, $F_4$\,,
$E_6$\,, $E_7$ and $E_8$ denote the compact connected simply connected groups
of those Cartan classification types, and their noncompact forms listed in
the Table are connected real forms contained as analytic subgroups 
in the corresponding
complex simply connected groups.  All 
known examples of compact connected simply connected quaternionic manifolds
with holonomy of quaternionic scalar type are riemannian symmetric spaces.

\begin{equation}\label{quaternion-symm}
\begin{tabular}{|l|l|c|c|} \hline
\multicolumn{4}{| c |}{Irreducible Quaternionic Symmetric Spaces, Scalar
        Part of Holonomy Quaternionic} \\ \hline \hline
\multicolumn{1}{|c}{compact $M = G/K$} &
        \multicolumn{1}{|c}{noncompact $M' = G'/K$} &
        \multicolumn{1}{|c}{Rank} &
        \multicolumn{1}{|c|}{Dimension/$\mathbb H$} \\ \hline \hline
$SU(r+2)/S(U(r)\times U(2))$ & $SU(r,2)/S(U(r)\times U(2))$ &
        $\min(r,2)$ & $r$ \\ \hline
$SO(r+4)/[SO(r) \times SO(4)]$ & $SO(r,4)/[SO(r) \times SO(4)]$ &
        $\min(r,4)$ & $r$  \\ \hline
$Sp(n+1)/[Sp(n) \times Sp(1)]$ & $Sp(n,1)/[Sp(n) \times Sp(1)]$ & $1$ & $n$
        \\ \hline
$G_2/SO(4)$ & $G_{2,A_1A_1}/SO(4)$ & $2$ & $2$ \\ \hline
$F_4/[Sp(3)\cdot Sp(1)]$ & $F_{4,C_3C_1}/[Sp(3)\cdot Sp(1)]$ & $4$ & $7$
        \\ \hline
$E_6/[SU(6) \cdot Sp(1)]$ & $E_{6,A_5C_1}/[SU(6) \cdot Sp(1)]$ & $4$ & $10$
        \\ \hline
$E_7/[Spin(12)\cdot Sp(1)]$ & $E_{7,D_6C_1}/[Spin(12)\cdot Sp(1)]$ &
        $4$ & $16$ \\ \hline
$E_8/[E_7\cdot Sp(1)]$ & $E_{8,E_7C_1}/[E_7\cdot Sp(1)]$ & $4$ & $28$ \\ \hline
\end{tabular}
\end{equation}
Thus irreducible quaternionic symmetric spaces have rank $1$, $2$, $3$ or $4$.
Curiously, quaternionic symmetric spaces for $F_4$\,, $E_6$\,, $E_7$\,, and
$E_8$ all have restricted root systems of type $F_4$\,.

\section{Complex Forms of Quaternionic Manifolds} \label{defs}
\setcounter{equation}{0}

Let $S$ be a smooth submanifold of a riemannian manifold $M$.
Let $A = \{A_x \mid x \in M\}$ denote a quaternionic structure on $M$.
If $x \in S$ let $A_x^S$ denote the subalgebra of all elements in
$A_x$ that preserve the real tangent space $T_x(S)$.
We say that $S$ is {\bf totally complex} if $A_x^S \cong \C$
and $T_x(S) \cap q(T_x(S)) = 0$ for all $q \in A_x \setminus A_x^s$\,,
for all $x \in S$.  If $S$ is totally complex in $M$ then
$A^S = \{A_x^S \mid x \in S\}$ restricts to a well defined almost complex
structure on $S$, parallel along $S$ because $A$ is parallel on $M$,
so (see \cite[Cor. 3.5, p. 145]{KN}) $(S, A^S|_S)$ is K\" ahler.
If in addition $\dim_\C S = \dim_\H M$ then we say that $S$ is
a {\bf maximal totally complex submanifold} of $M$.
\smallskip

Let $S$ be a maximal totally complex submanifold of $M$.  Suppose that $S$ 
is a topological component of the fixed point set of an involutive
isometry $\sigma$ of $M$.  Then we say that $S$ is a {\bf complex form}
of $M$ and that $\sigma$ is the {\bf quaternion conjugation} of $M$
over $S$.  The following is immediate.

\begin{lemma} \label{lem1}
Let $(M,A)$ be a quaternionic symmetric space.  If $S$ is a complex form
of $M$, then $S$ is a totally geodesic submanifold.  If $S$ is a totally
geodesic, totally complex submanifold of $M$, then $S$ is an hermitian
symmetric space.
\end{lemma}

Let $M = G/K$, irreducible quaternionic symmetric space, with base point 
$x_0 = 1K$, where $K = K'\cdot Sp(1)$ as in Proposition 
\ref{quaternion-criterion}(iii) and Table \ref{quaternion-symm}.
Let $\theta$ denote the involutive automorphism of $G$ that is conjugation
by the symmetry (say $t$) at $x_0$\,.
Let $S \subset M$ be a totally geodesic submanifold
through $x_0$\,.  Then $S$ is a riemannian symmetric space with symmetry
$t|_{_S}$ at $x_0$\,. Express $S = L(x_0) \cong L/V$ where 
$L$ is the identity component of $\{g \in G \mid g(S) = S\}$ and 
$V = L \cap K$.  Then $\theta(L) = L$.
\medskip

The following three results are our basic tools for finding the complex forms
$S = L/V$ of $M = G/K$ where rank$(L)$ = rank$(G)$.  Proposition 
\ref{equal-rank-criterion} gives criteria for $L/V$ to be an appropriate
submanifold of $G/K$, Proposition \ref{equal-rank-exist} tells us that when
$L/V$ is identified abstractly it in fact exists well positioned in $G/K$, and
Proposition \ref{equal-rank-distinguish} is a uniqueness theorem
showing just when two complex forms are $G$--equivalent.

\begin{proposition} \label{equal-rank-criterion}
Let $M = G/K$ be an irreducible quaternionic symmetric space, with base point
$x_0 = 1K$, as above.  Let $\sigma$ be an involutive inner automorphism of $G$
that commutes with $\theta$.  Let $L$ be the identity component of the fixed
point set $G^\sigma$.  Set $V = L \cap K$.  Denote $S = L(x_0) \cong L/V$.

{\rm 1.} If $V \cap Sp(1)$ is a circle group then $S$ is a totally complex
submanifold of $M$.

{\rm 2.} $S$ is a complex form of $M$ if and only if {\rm (i)} $V \cap Sp(1)$
is a circle group and {\rm (ii)} $\dim_\C S = \dim_\H M$.

{\rm 3.} If $S$ is a complex form of $M$ then $\sigma = \text{ Ad}(s)$ where
$s \in V$.
\end{proposition}

\begin{proposition} \label{equal-rank-exist}
Let $M = G/K$ be an irreducible quaternionic symmetric space, with base point
$x_0 = 1K$, as above.  Let $L$ be a symmetric subgroup of equal rank in $G$
that has an hermitian symmetric quotient $L/V$ such that $V$ is isomorphic
to a symmetric subgroup $V' \subset K$.  Then $L$ is conjugate to a
$\theta$--stable subgroup $L' \subset G$ such that $L' \cap K = V'$.
\end{proposition}

\begin{proposition} \label{equal-rank-distinguish}
Let $M = G/K$ be an irreducible quaternionic symmetric space, with base point
$x_0 = 1K$, as above.  Let $S_i = L_i(x_0) \cong K_i/V_i$ be two complex forms
of $M$.  If $S_1$ and $S_2$ are isometric then some element of $K$
carries $S_1$ onto $S_2$\,.
\end{proposition}

\noindent {\bf Proof of Proposition \ref{equal-rank-criterion}.}
We can pass to the compact dual if necessary, so we may (and do) assume
$M$ compact.

Decompose the Lie algebra $\gg$ of $G$ under $d\theta$, $\gg = \gk + \gm$
where $\gk$ is the Lie algebra of $K$ and $\gm$ represents the real tangent
space of $M$.  Then $Sp(1)$ gives $\gm$ a quaternionic vector space structure,
so any circle subgroup gives $\gm$ a complex vector space structure.  If
that circle is $V \cap Sp(1)$ it defines an $L$--invariant almost complex
structure on $S$, and that is integrable because $S$ is a riemannian symmetric
space.  We have proved Statement 1.
\smallskip

For Statement 2 first suppose that $S$ is a complex form of $M$.  Since
$\sigma$ is inner by hypothesis, rank$(L)$  = rank$(G)$.  Since $S$ is an
hermitian symmetric space, rank$(V)$ = rank$(L)$.  Now $V$ contains a
Cartan subgroup $T$ of $G$.  Thus $V \cap Sp(1)$ contains a circle group
$T_1 := T \cap Sp(1)$.  Now the only possibilities for $V \cap Sp(1)$ are
(a) $T_1$\,, (b) the normalizer of $T_1$ in $Sp(1)$, and (c) all of $Sp(1)$.
Here (b) is excluded because it would prevent $S$ from having an $L$--invariant
almost complex structure, and (c) is excluded because it would prevent $S$
from being totally complex, so $V \cap Sp(1)$ is a circle group.  Finally
$\dim_\C S = \dim_\H M$ because $S$ is a maximal totally complex submanifold
of $M$.
\smallskip

Conversely suppose that $V \cap Sp(1)$ is a circle group and 
$\dim_\C S = \dim_\H M$.  By Statement 1, $S$ is a totally complex
submanifold of $M$.  By $\dim_\C S = \dim_\H M$ it is a maximal totally
complex submanifold.  And we started with the symmetry $\sigma$, so
$S$ is a complex form of $M$.  
\smallskip

For Statement 3 note, as above, that $s \in L$ because rank$(L)$  = rank$(G)$,
and now $s \in V$ because rank$(V)$ = rank$(L)$.
\hfill $\square$
\medskip

\noindent {\bf Proof of Proposition \ref{equal-rank-exist}.}
All our groups have equal rank, so $V'$ is the $K$--centralizer of some
$v' \in V'$ with $v'^2$ central in $K$.  Here $K$ contains the center of $G$,
and those centers satisfy $Z_K/Z_G = \{1,z\}Z_G$ cyclic order $2$.
Let $\sigma' = \text{ Ad}(v')$.  If $v'^2 \in zZ_G$ then $\sigma'^2 = \theta$,
so $d\sigma$ has eigenvalues $\pm \sqrt{-1}$ on $\gm$, and $L' = G^{\sigma'}$
has the property that $S' = L'(x_0) \cong L'/V'$ is hermitian symmetric.
Since $V \in L$ and $V' \in L'$ are symmetric subgroups of $G$, and their
hermitian symmetric subgroups are isomorphic, it follows from Table 
\ref{quaternion-symm} and the classification of riemannian symmetric spaces 
that $L \cong L'$\,.  Now $L$ and $L'$ are conjugate in $G$ so we may assume
$L = L'$.  then $V$ and $V'$ are isomorphic symmetric subgroups in $L$, so
they are $L$--conjugate.  This completes the proof.  \hfill $\square$

\medskip

\noindent {\bf Proof of Proposition \ref{equal-rank-distinguish}.}
Suppose that $S_1$ and $S_2$ are isometric, say $g: S_1 \cong S_2$ for
some isometric map $g$.  We can assume $g(x_0) = x_0$\,,
so $dg$ gives a Lie triple system isomorphism of $\gl_1 \cap \gm$
onto $\gl_2 \cap \gm$.  Write $\gl_i = \gl'_i \oplus \gz_i$ where
$\gl'_i$ is generated by $\gl_i \cap \gm$ and $\gz_i \subset \gv_i$ is
a complementary ideal.  Then $dg$ gives a Lie algebra isomorphism of 
$\gl'_1$ onto $\gl'_2$\,.  Let ${\mathbf j}_i \in \gs \gp (1)$ be orthogonal
to the Lie algebra of the circle group $V_i \cap Sp(1)$.  Then 
${\mathbf j}_i$ centralizes $\gz_i$ and $\gm$ is the real vector space 
direct sum of $\gl_i \cap \gm$ with ad$({\mathbf j}_i)(\gl_i \cap \gm)$.
Now ad$(\gz_i)|_\gm = 0$, so each $\gz_i = 0$, and $dg : \gl_1 \cong \gl_2$\,.
Since $\gl_1$ and $\gl_2$ are isomorphic symmetric subalgebras of $\gg$ they
are Ad$(G)$--conjugate.  Thus we may assume $g \in G$.  As $g(x_0) = x_0$
now $g \in K$.  Thus some $g \in K$ carries $S_1$ onto $S_2$\,.
\hfill $\square$
\medskip

Propositions \ref{equal-rank-criterion} and \ref{equal-rank-distinguish} 
will let us
do the classification of complex forms $S = L/V$ of quaternionic symmetric
spaces $M = G/K$ in case rank$(L)$ = rank$(G)$.  There are only a few
cases where rank$(L) <$ rank$(G)$, and we will handle them individually.
That is not very elegant, but it is very efficient.

\section{The Classification of Complex Forms.} \label{classification}
\setcounter{equation}{0}

In this section we state the classification of complex forms $S = L/V$
of quaternionic symmetric spaces $M = G/K$ and $M' = G'/K$ whose holonomy
has quaternion scalar part.  The proofs are given in Sections
\ref{equal-rank-classical}, \ref{equal-rank-exceptional} and
\ref{unequal-rank}.  We state the results separately for the compact and
the noncompact cases.

\begin{theorem} \label{class-cplx-forms-compact}
Let $M = G/K$ be a compact simply connected irreducible quaternionic 
riemannian symmetric space.  Then the complex forms $S = L/V$ of $M$
are exactly the following, and each is unique up to the action
of $G$.
\medskip

\noindent {\rm \bf 1.} $M = SU(r+2)/S(U(r)\times U(2))$.  Then 
{\rm (1a)} $S = SO(r+2)/[SO(r) \times SO(2)]$, or \hfill\newline
{\rm (1b)} $S =
[SU(u+1)/S(U(u) \times U(1))] \times [SU(r-u+1)/S(U(r-u) \times U(1))]$,
$0 \leqq u \leqq r$.
\medskip

\noindent {\rm \bf 2.} $M = SO(r+4)/[SO(r) \times SO(4)]$.  Then 
{\rm (2a)} 
$S = SU(r'+2)/S(U(r') \times U(2))$, $r = 2r'$ even, or
{\rm (2b)}
$S = \{SO(u+2)/[SO(u) \times SO(2)]\} \times 
\{SO(r-u+2)/[SO(r-u) \times SO(2)]\}$, $0 \leqq u \leqq r$.
\medskip

\noindent {\rm \bf 3.} $M = Sp(n+1)/[Sp(n) \times Sp(1)] = P^n(\H)$.  Then
$S = U(n+1)/[U(n) \times U(1)] = P^n(\C)$.
\medskip

\noindent {\rm \bf 4.} $M = G_2/SO(4)$.  Then $S = SO(4)/[SO(2) \times SO(2)]
=  P^1(\C) \times  P^1(\C)$.
\medskip

\noindent {\rm \bf 5.} $M = F_4/[Sp(3)\cdot Sp(1)]$.  Then $S = 
\{Sp(3)/U(3)\} \times P^1(\C)$.
\medskip

\noindent {\rm \bf 6.} $M = E_6/[SU(6)\cdot Sp(1)]$.  Then 
{\rm (6a)} $S = SO(10)/U(5)$, or
{\rm (6b)} 
$S = Sp(4)/U(4)$, or
\hfill\newline 
{\rm (6c)} 
$S = [SU(6)/S(U(3) \times U(3))] \times P^1(\C)$.
\medskip

\noindent {\rm \bf 7.} $M = E_7/[Spin(12)\cdot Sp(1)]$.  Then 
{\rm (7a)} $S = E_6/[Spin(10)\cdot U(1)]$, or
\hfill\newline
{\rm (7b)} $S = SU(8)/S(U(4) \times U(4))$, or
{\rm (7c)} $S = SO(12)/U(6) \times P^1(\C)$.
\medskip

\noindent {\rm \bf 8.}  $M = E_8/[E_7\cdot Sp(1)]$.  Then 
{\rm (8a)} $S = [E_7/E_6T_1] \times P^1(\C)$ or
{\rm (8b)} $S = SO(16)/U(8)$.
\end{theorem}

\begin{theorem} \label{class-cplx-forms-noncompact}
Let $M = G/K$ be a noncompact irreducible quaternionic
riemannian symmetric space.  Then the complex forms $S = L/V$ of $M$
are exactly the following, and each is unique up to the action
of $G$.
\medskip
 
\noindent {\rm \bf 1.} $M = SU(r,2)/S(U(r)\times U(2))$.  Then
{\rm (1a)} $S = SO(r,2)/[SO(r) \times SO(2)]$, or \hfill\newline
{\rm (1b)} $S =
[SU(u,1)/S(U(u) \times U(1))] \times [SU(r-u,1)/S(U(r-u) \times U(1))]$,
$0 \leqq u \leqq r$.
\medskip
 
\noindent {\rm \bf 2.} $M = SO(r,4)/[SO(r) \times SO(4)]$.  Then
{\rm (2a)}
$S = SU(r',2)/S(U(r') \times U(2))$, $r = 2r'$ even, or
{\rm (2b)}
$S = \{SO(u,2)/[SO(u) \times SO(2)]\} \times
\{SO(r-u,2)/[SO(r-u) \times SO(2)]\}$, $0 \leqq u \leqq r$.
\medskip
 
\noindent {\rm \bf 3.} $M = Sp(n,1)/[Sp(n) \times Sp(1)] = H^n(\H)$.  Then
$S = U(n,1)/[U(n) \times U(1)] = H^n(\C)$.
\medskip
 
\noindent {\rm \bf 4.} $M = G_{2,A_1A_1}/SO(4)$.  Then 
$S = SO(2,2)/[SO(2) \times SO(2)] = H^1(\C) \times H^1(\C)$.
\medskip
 
\noindent {\rm \bf 5.} $M = F_{4,C_3C_1}/[Sp(3)\cdot Sp(1)]$.  Then $S =
\{Sp(3;\R)/U(3)\} \times H^1(\C)$.
\medskip
 
\noindent {\rm \bf 6.} $M = E_{6,A_5C_1}/[SU(6)\cdot Sp(1)]$.  Then
{\rm (6a)} $S = SO^*(10)/U(5)$, 
{\rm (6b)}
$S = Sp(4;\R)/U(4)$, or
\hfill\newline
{\rm (6c)}
$S = [SU(3,3)/S(U(3) \times U(3))] \times H^1(\C)$.
\medskip
 
\noindent {\rm \bf 7.} $M = E_{7,D_6C_1}/[Spin(12)\cdot Sp(1)]$.  Then
{\rm (7a)} $S = E_{6,D_5T_1}/[Spin(10)\cdot U(1)]$, or
\hfill\newline
{\rm (7b)} $S = SU(4,4)/S(U(4) \times U(4))$, or
{\rm (7c)} $S = SO^*(12)/U(6) \times H^1(\C)$.
\medskip
 
\noindent {\rm \bf 8.}  $M = E_{8,E_7C_1}/[E_7\cdot Sp(1)]$.  Then
{\rm (8a)} $S = [E_{7,E_6T_1}/E_6T_1] \times H^1(\C)$ or
{\rm (8b)} $S = SO^*(16)/U(8)$,
\end{theorem}

Of course Theorem \ref{class-cplx-forms-noncompact} is immediate from
Theorem \ref{class-cplx-forms-compact} by passage to the compact dual 
symmetric spaces, so we need only prove Theorem \ref{class-cplx-forms-compact}.
The proof of Theorem \ref{class-cplx-forms-compact} consists of
consolidating the results of Sections \ref{equal-rank-classical}, 
\ref{equal-rank-exceptional} and \ref{unequal-rank}.

\section{The Equal Rank Classification --- Classical Cases} 
\label{equal-rank-classical}
\setcounter{equation}{0}

We run through the list of compact irreducible quaternionic symmetric spaces
$M = G/K$ from Table \ref{quaternion-symm}, for the cases where $G$ is a
classical group.  For each of them we look at
the possible symmetric subgroups $L$ that correspond to an hermitian symmetric
space $S = L/V$ such that rank$(L)$ = rank$(G)$, $\dim_\C S = \dim_\H M$, 
rank$(S) \leqq$ rank$(M)$, and $V$ is isomorphic to a symmetric subgroup
of $K$ properly places as in Proposition \ref{equal-rank-criterion}.  The equal 
rank classification will follow using Proposition \ref{equal-rank-distinguish}.
We retain the notation used in Propositions \ref{equal-rank-criterion} and 
\ref{equal-rank-distinguish}.  Fix $s \in K$ such that $L$ is the identity
component of $\sigma = \text{ Ad}(s)$.
\medskip

{\sc Case} $M = SU(r+2)/S(U(r) \times U(2))$.  First suppose $r \geqq 2$.  
We may take $s$ to be diagonal.  It has only two distinct eigenvalues, and
its component in the $U(2)$--factor of $K$ must have both eigenvalues.
Now $L \cong S(U(u+1) \times U(v+1))$, 
$V \cong S( [U(u) \times U(1)]  \times [U(v) \times U(1)])$, and
$S$ is the product $P^u(\C) \times P^v(C)$ of complex projective spaces.
Here $\dim_\H M = r = u+v = \dim_\C S$. If $u, v \geqq 1$ then
rank$(M) = 2 =$ rank$(S)$.  If $u = 0$ then the factor $P^u(\C)$ is reduced to
a point, $S \cong P^v(C)$, and rank$(S) = 1$.  The analog holds, of course, if
$v = 0$.
\smallskip

Now consider the degenerate case $r = 1$.  Then $M = P^2(\C)$ and fits the
dual pattern described in the paragraph just before the statement of 
Proposition \ref{quaternion-criterion}.  Relative to the quaternionic
structure denoted $A_1$ there, the one with with quaternion scalar part,
the matrix considerations above show that $M$ has a complex form
$S = P^1(\C)$.
\medskip

{\sc Case} $M = SO(r+4)/[SO(r) \times SO(4)]$.  
As before, the matrix $s$ has just two distinct eigenvalues, and both
each must appear with multiplicity $2$ in the $SO(4)$--factor of $K$.  
If $s^2 = I$ then $L = SO(u+2) \times SO(v+2)$ with $u+v = r$,
where $V = L \cap K = [SO(u) \times SO(2)] \times [SO(v) \times SO(2)]$.
here the $SO(2)$--factors in $V$ are the intersection with the $SO(4)$--factor
of $K$.  That gives us the complex forms
$S = (SO(u+2)/[SO(u) \times SO(2)]) \times (SO(v+2)/[SO(v) \times SO(2)])$
of $M$.
\smallskip

Now suppose $s^2 = -I$.  Then
$r = 2r'$  even, $L \cong U(r'+2)$, $V \cong U(r') \times U(2)$, and 
we have the complex form
$S \cong SU(r'+2)/S(U(r') \times U(2))$ of $M$.
\medskip

{\sc Case} $M = Sp(n+1)/[Sp(n) \times Sp(1)] = P^n(\H)$.  The symmetric 
subgroups of
$Sp(n+1)$ are the $Sp(u) \times Sp(v)$, $u+v = n+1$, and $U(n+1)$.  The
first case, $L = Sp(u) \times Sp(v)$, would give 
$V = Sp(u) \times Sp(v-1) \times Sp(1)$, so $S = Sp(v)/[Sp(v-1) \times Sp(1)]$,
which is not hermitian symmetric.  That leaves the case $L = U(n+1)$ and
$V = U(n) \times U(1)$ where $S = P^n(\C)$.  It satisfies the conditions of
Proposition \ref{equal-rank-criterion} and thus is a complex form of $M$.

\section{The LiE Program} \label{prog}
\setcounter{equation}{0}

While the matrix computation methods of Section \ref{equal-rank-classical}
work well for the classical group cases, it is more convenient to make use
of the root structure in the exceptional group cases.  In this section we
indicate just how we used the LiE program \cite{L} to do that.  We illustrate
it for $E_8$\,, but it is the same for any simple group structure.
Here node refers to the simple root at which the negative of the maximal
root is attached in the extended Dynkin diagram.
\medskip

\noindent {\sc Step 0: Initialize.}
\smallskip

$>$ setdefault(E8)

$>$ rank = 8

$>$ diagram  \phantom{XXXX} ; {\em prints out the Dynkin diagram and numbers the simple roots.}

$>$ node = 8 \phantom{XXXi} ; {\em the number of the simple root that defines $K$.}
\medskip

\noindent {\sc Step 1: Positive Roots of $\gg$.}
\smallskip

$>$ pos = pos\_roots

$>$ max\_root = pos[n\_rows(pos)]
\medskip

\noindent {\sc Step 2: Positive Roots of $\gk$.}
\smallskip

$>$ kkk = pos
 
$>$ for i = 1 to n\_rows(kkk) do \\
\phantom{XXX} if kkk[i,node] == 1 then kkk[i] = null(rank) fi od
\phantom{XXX} ; {\em zeroes rows for roots of $\gm$}
 
$>$  kk = unique(kkk)
\phantom{XXXXXXXXXXXXXXXXXXXXXi} ; {\em eliminates duplicate rows}
 
$>$  k = null(n\_rows(kk)-1,rank)
 
$>$ for i = 1 to n\_rows(k) do k[i] = kk[i+1] od
\phantom{XXXXXXX} ; {\em eliminates last zero row}
 
$>$ Cartan\_type(k)
\phantom{XXXXXXX} 
; {\em verifies correct Cartan type for $\gk$, in this case $E_7A_1$}
\medskip

\noindent {\sc Step 3:  Positive Roots of $\gm$.}
\smallskip

$>$ mmm = pos
 
$>$ for i = 1 to n\_rows(mmm) do \\
\phantom{XXX}  if mmm[i,node] != 1 then mmm[i] = null(rank) fi od
\phantom{XXX} ; {\em zeroes rows for roots of $\gk$}
 
$>$ mm = unique(mmm)
\phantom{XXXXXXXXXXXXXXXXXXXXX} ; {\em eliminates duplicate rows}
 
$>$ m = null(n\_rows(mm)-1,rank)
 
$>$ for i = 1 to n\_rows(m) do m[i] = mm[i+1] od
\phantom{XXXXXXX} ; {\em eliminates last zero row}
\medskip
 
\noindent {\sc Step 4:  Choice of $sym$ where $\sigma = \text{\rm Ad}(sym)$;
definition of $\gl = \gg^\sigma$.}
\smallskip

$>$ sym = null(rank + 1) \phantom{XXXX} ; {\em initializes $sym$ as row vector}

$>$ sym[node] = 1 \phantom{XXXXXXXXi} 
; {\em one possibility for nonzero element of $sym$}

$>$ sym[rank+1] = 2 \phantom{XXXXXXX} 
; {\em normalizes $1$--parameter group containing symmetry $sym$}
 
$>$ l = cent\_roots(sym) \phantom{XXXXXi} 
; {\em defines $\gl$ as centralizer of $sym$}
 
$>$ Cartan\_type(l) \phantom{XXXXXXXW} 
; {\em Cartan type of $\gl$, in this case $E_7A_1$}
\medskip

\noindent {\sc Step 5:  Positive Roots of $\gs := \gl \cap \gm$ and 
of $\gv := \gl \cap \gk$.}
\smallskip

$>$ sss = l
 
$>$ for i = 1 to n\_rows(sss) do 
  if sss[i,node] != 1 then sss[i] = null(rank) fi od
 
$>$ ss = unique(sss)
 
$>$ s = null(n\_rows(ss)-1,rank)
 
$>$ for i = 1 to n\_rows(s) do s[i] = ss[i+1] od
 
$>$ vvv = l
 
$>$ for i = 1 to n\_rows(vvv) do 
  if vvv[i,2] == 1 then vvv[i] = null(rank) fi od
 
$>$ vv = unique(vvv)
 
$>$ v = null(n\_rows(vv)-1,rank)
 
$>$ for i = 1 to n\_rows(v) do v[i] = vv[i+1] od
 
$>$ Cartan\_type(v) \phantom{XXXXXXi}
; {\em Cartan type of $\gv$, in this case $E_6T_1T_1$}
\medskip

\noindent ; {\em At this point we have the information that $S = L/V
\cong (E_7/[E_6 \times T_1]) \times (T_1 / T_1)$}, \\
; {\em so it is an hermitian symmetric subspace of $G/K$.}
\medskip

\noindent {\sc Step 6:  Verify that $S$ is a maximal totally complex 
submanifold of $M$.}
\smallskip

$>$ t =  null(n\_rows(s)-1,rank) 
 
$>$ for i=1 to  n\_rows(t) do t[i] = max\_root - s[i] od
 
$>$ u = null(n\_rows(s) + n\_rows(t) + n\_rows(m), rank)
 
$>$ for i = 1 to n\_rows(s) do u[i] = s[i] od
 
$>$ for i = 1 to n\_rows(t) do u[n\_rows(s) + i] = t[i] od
 
$>$ for i = 1 to n\_rows(m) do 
  u[n\_rows(s) + n\_rows(t) + i] = m[i] od \\
\\
\phantom{XX}; {\em now the rows of u are: positive roots of $\gs$}, \\
\phantom{XX}; \phantom{XXXXXXXXXXXXX} {\em maximal root minus positive 
roots of $\gs$}, \\ 
\phantom{XX}; \phantom{XXXXXXXXXXXXX} {\em positive roots of $\gm$}
\medskip
 
$>$ w = unique(u) 
\phantom{XXi} ; {\em the rows of w are the positive roots of $\gm$ and 
non--root linear} \\ 
\phantom{XXXXXXXXXXXXXXX} ; {\em functionals {\rm (} max root 
minus positive root of $\gs$ {\rm )}}
\\
 
$>$ n\_rows(w) - n\_rows(m)
\phantom{XX} ; {\em number of non--root linear functionals in w}, \\
\phantom{XXXXXXXXXXXXXXXXXXXi} ; {\em measures failure of $S$ to be maximal totally complex};\\
\phantom{XXXXXXXXXXXXXXXXXXXi} ; {\em OK here because it returns $0$}
\bigskip

We carry out the routine in some key cases.  These are cases 
where $K$ and $L$ are conjugate in $G$.  
\medskip

{\sc Case} $G = B_7$\,.  Here node $= 2$, and sym = $[0,1,0,0,0,0,0,2]$ leads
to $L = B_5A_1A_1$ and $V = B_4T_1T_1T_1$\,, thus to the complex form
$S = SO(11)/[SO(9) \times SO(2)] \times  P^1(\C) \times P^1(\C)$
of $G/K = SO(15)/[SO(11) \times SO(4)$.  More generally, for $B_n$ with
$n \geqq 3$,  node $= 2$, and sym = $[0,1, 0, \dots , 0, 2]$ gives the complex
form $S = SO(2n-3)/[SO(2n-5) \times SO(4)] \times  P^1(\C) \times P^1(\C)$
of $G/K = SO(2n+1)/[SO(2n-3) \times SO(4)]$.  This is the case $v = 2,
u = r-2$ considered for $G =  SO(r+4)$, $r$ odd, 
in Section \ref{equal-rank-classical}.
\medskip

{\sc Case} $G = D_7$\,.  Here node $= 2$, and sym = $[0,1,0,0,0,0,0,2]$ leads
to $L = D_5A_1A_1$ and $V = D_4T_1T_1T_1$\,, thus to the complex form
$S = SO(10)/[SO(8) \times SO(2)] \times  P^1(\C) \times P^1(\C)$
of $G/K = SO(14)/[SO(10) \times SO(4)]$.  More generally, for $D_n$ with
$n \geqq 3$,  node $= 2$, and sym = $[0,1, 0, \dots , 0, 2]$ gives the complex
form $S = SO(2n-4)/[SO(2n-6) \times SO(2)] \times  P^1(\C) \times P^1(\C)$
of $G/K = SO(2n)/[SO(2n-4) \times SO(4)]$.  This is the case $v = 2,
u = r-2$ considered for $G =  SO(r+4)$, $r$ even, 
in Section \ref{equal-rank-classical}.
\medskip

{\sc Case} $G = G_2$\,.  Here node $= 2$, and sym = $[0,1,2]$ leads to
$L = A_1A_1$ and $V = T_1T_1$\,, thus to the complex form
$S = P^1(\C) \times P^1(\C)$ of $G/K = G_2/SO(4)$.
\medskip

{\sc Case} $G = F_4$\,.  Here node $= 1$, and sym = $[1,0,0,0,2]$ leads to
$L = C_3C_1$ and $V = A_2T_1T_1$\,, thus to the complex form
$S = [Sp(3)/U(3)] \times  P^1(\C)$ of $G/K = F_4/C_3C_1$\,.
\medskip

{\sc Case} $G = E_6$\,.  Here node $= 2$, and sym = $[0,1,0,0,0,0,2]$ leads
to $L = A_5A_1$ and $V = A_2T_1A_2T_1$\,,  thus to the complex form
$S = [SU(6)/S(U(3) \times U(3))] \times  P^1(\C)$ of $G/K = E_6/A_5A_1$\,.
\vskip .15 cm

{\sc Case} $G = E_7$\,.  Here node $= 2$, and sym = $[0,1,0,0,0,0,0,2]$ leads
to $L = D_6A_1$ and $V = A_5T_1T_1$\,, and thus to the complex form 
$S = [SO(12)/U(6)] \times P^1(\C)$ of $G/K = E_7/D_6A_1$\,.
\medskip

{\sc Case} $G = E_8$\,.  As we saw, sym = $[0,0,0,0,0,0,0,1,2]$ leads to
$L = E_7A_1$ and $V = E_6T_1T_1$\,, and thus to the complex form 
$S = (E_7/[E_6 \times T_1]) \times P^1(\C)$ of $G/K = E_8/E_7A_1$\,.
\medskip

{\sc Cases} $A_7$ and $C_7$\,.  Here the computation using LiE has not
yet produced complex forms $S$ of $M$.  In other words
I have not yet guessed the appropriate vectors sym to define toral elements
of $G$ whose centralizers are appropriate subgroups $L \subset G$.

\section{The Equal Rank Classification --- Exceptional Cases} 
\label{equal-rank-exceptional}
\setcounter{equation}{0}

In this section we complete the classification for the equal rank 
exceptional group cases.  

{\sc Case} $G = G_2$\,.  The only symmetric subgroup of $G_2$ is $SO(4)$,
so here the only complex form of $M = G_2/SO(4)$ is 
$S =  P^1(\C) \times P^1(\C)$ as described in Section \ref{prog}.
\medskip

{\sc Case} $G = F_4$\,.  The only symmetric subgroups of $F_4$ are
$Sp(3)\cdot Sp(1)$ and $Spin(9)$.  If $L = Spin(9)$ then the hermitian
symmetric space $L/V = Spin(9)/[Spin(7) \times Spin(2)]$.  That would
place the $Spin(7)$--factor of $V$ in the $Sp(3)$--factor of $K$; but
$Sp(3) \subset SU(6)$ while $Spin(7)$ has no nontrivial representation of
degree $< 7$.  Thus $L \ne Spin(9)$, so here the only complex
form of $M = F_4/C_3C_1$ is $S = [Sp(3)/U(3)] \times  P^1(\C)$ as
described in Section \ref{prog}.
\medskip

{\sc Case} $G = E_6$\,.  The only symmetric subgroups of maximal rank in 
$E_6$ are $A_5A_1$ and $D_5T_1$\,.  
\smallskip

If $L = D_5T_1$
then the hermitian symmetric space $S = L/V$ must be
$SO(10)/[SO(8) \times SO(2)]$ with $V = [SO(8) \times SO(2)]\cdot SO(2)$,
or $[SO(10)/U(5)]$ with $V = U(5)\cdot SO(2)$.  The first of these is 
excluded because $\dim_\C SO(10)/[SO(8) \times SO(2)] = 8 < 10 = \dim_\H M$.
The second of these is a complex form of $M = E_6/A_5A_1$ by Propositions
\ref{equal-rank-criterion} and \ref{equal-rank-exist}.
\smallskip

$L = A_5A_1$ gives another complex form 
$S = [SU(6)/S(U(3) \times U(3))] \times  P^1(\C)$ of
$M = E_6/A_5A_1$ as described in Section \ref{prog}.
\medskip

{\sc Case} $G = E_7$\,.  The only symmetric subgroups of $E_7$ are
$D_6A_1$\,, $A_7$ and $E_6T_1$\,.  

\smallskip
If $L \cong E_6T_1$ then
the hermitian symmetric space $S = L/V$ must be $E_6/D_5T_1$ with 
$V = D_5T_1T_1$\,.  It is a complex form of $M =  E_7/D_6A_1$ by
Propositions \ref{equal-rank-criterion} and \ref{equal-rank-exist}.
\smallskip

If $L = A_7$ then the hermitian symmetric space $S = L/V$ must be 
$SU(8)/S(U(u) \times U(v))$ with $u + v = 8$.  Here $\dim_\C L/V = UV$ 
while $\dim_\H M = 16$, so $u = v = 4$.  That would place the 
$[SU(4) \times SU(4)]$--factor of $V$ in the $Spin(12)$--factor of $K$.
It could only sit there as $Spin(6) \times Spin(6)$, which is the identity
component of its $Spin(12)$--normalizer because it is a symmetric
subgroup of $Spin(12)$, so the circle center of $V$ is contained in the
$Sp(1)$--factor of $K$.  Thus $S$ is a complex form of $M =  E_7/D_6A_1$ by
Propositions \ref{equal-rank-criterion} and \ref{equal-rank-exist}.
\smallskip

$L = D_6A_1$  gives another complex form $S = [SO(12)/U(6)] \times P^1(\C)$ 
of $M =  E_7/D_6A_1$ as as described in Section \ref{prog}.
\medskip

{\sc Case} $G = E_8$\,.  The only symmetric subgroups of $E_8$ are
$E_7A_1$ and $D_8$.  
\smallskip

If $L = D_8$ then the hermitian symmetric 
space $S = L/V$ either must be $SO(16)/[SO(14) \times SO(2)]$ with
$V = [SO(14) \times SO(2)]$, or $SO(16)/U(8)$ with $V = U(8)$.
The first of these is excluded because 
$\dim_\C SO(16)/[SO(14) \times SO(2)] = 14 < 28 = \dim_\H M$.
The second of these is a complex form of $M = E_8/E_7A_1$ by
Propositions \ref{equal-rank-criterion} and \ref{equal-rank-exist}.
\smallskip

$L \cong E_7A_1$ gives another complex form
$S = (E_7/[E_6 \times T_1]) \times P^1(\C)$ of $M = E_8/E_7A_1$
as described in Section \ref{prog}.

\section{The Unequal Rank Classification} \label{unequal-rank}
\setcounter{equation}{0}

In this section we deal with the cases rank$(L) <$ rank$(G)$.
Here $G$ is of type $A_n$\,, $D_n$ or $E_6$\,.
\medskip

{\sc Case} $M = SU(r+2)/S(U(r) \times U(2))$.  The only symmetric 
subgroups of lower rank in $SU(r+2)$ are $SO(r+2)$ and, for $r = 2r'$ even,
$Sp(r'+1)$.  
\smallskip

If $L = Sp(r'+1)$, $r = 2r'$ even, then the hermitian symmetric space 
$S = Sp(r'+1)/U(r'+1)$ with $V = U(r'+1)$.  Here $\dim_\H M = 2r'$ and 
$\dim_\C S = \frac{1}{2}(r'+2)(r'+1)$, so those dimensions are equal just
when $r'^2 - r' + 2 = 0$.  That equation has no integral solution.
Thus $L \ne Sp(r'+1)$.
\smallskip

If $L = SO(r+2)$ then the hermitian symmetric space 
$S = SO(r+2)/[SO(r) \times SO(2)]$ with $V = [SO(r) \times SO(2)]$.
The $SO(2)$--factor of $V$ is contained in the derived group $SU(2)$
of the $U(2)$--factor of $K$, and $\dim_\C S = r =  \dim_\H M$, 
so Proposition \ref{equal-rank-criterion} shows that $S$ is a complex form
of $M$.
\medskip

{\sc Case} $M = SO(2n+4)/[SO(2n) \times SO(4)]$.  The only symmetric
subgroups of lower rank in $SO(2n+4)$ are $SO(2u+1) \times SO(2v+1)$
where $u + v = n+1$.  If $L = SO(2u+1) \times SO(2v+1)$ then $S =
\{SO(2u+1)/[SO(2u-1)\times SO(2)]\}\times \{SO(2v+1)/[SO(2v-1)\times SO(2)]\}$
with $V = SO(2u-1)\times SO(2)\times SO(2v-1)\times SO(2)$, where the
product of the two $SO(2)$--factors is contained in the $SO(4)$--factor of
$K$.  Since $\dim_\C S = (2u-1)+(2v-1) = n = \dim_\H M$, the argument of
Proposition \ref{equal-rank-criterion} shows that $S$ is a complex form
of $M$.
\medskip

{\sc Case} $M = E_6/A_5A_1$\,.  The only symmetric subgroups of lower rank
in $E_6$ are $F_4$ and $C_4$\,, and $L \ne F_4$ because $F_4$ has no hermitian
symmetric quotient space.  If $L = Sp(4)$ then $S = Sp(4)/U(4)$ with $V = U(4).$
Here $V$ sits in $K$ as follows.  The semisimple part $[V,V] = U(4)/\{\pm I\}
= SO(6) \subset SU(6) = A_5$\,.  $[V,V]$ is a connected symmetric subgroup 
of the  connected simple group $A_5$\,, so it is equal to the identity 
component of its normalizer in $A_5$\,.  Thus the projection
$K = A_5A_1 \to A_5$ annihilates the circle center of $V$.  In other words,
$V \cap Sp(1)$ is a circle group central in $V$.  It follows as in
Proposition \ref{equal-rank-criterion}(1) that $S$ is a totally complex
submanifold of $M$.  Since $\dim_\C S = 10 = \dim_\H M$ it is a maximal
totally complex submanifold, and being a symmetric submanifold it is a
complex form.
\bigskip

This completes the proof of Theorems \ref{class-cplx-forms-compact}
and \ref{class-cplx-forms-noncompact}, the main results of this note.

\section{Quaternionic Forms} \label{quat-form}
\setcounter{equation}{0}

In this section we look at the idea of quaternionic forms of 
symmetric spaces as suggested by the examples
of projective planes $P^2(\H) \subset P^2({\mathbb O})$
and hyperbolic planes $H^2(\H) \subset H^2({\mathbb O})$.  
The meaning of Cayley structure is not entirely clear
because of nonassociativity, 
so we do not have a good definition for Cayley symmetric space.
Here we offer a tentative definition of quaternionic form and a number
of examples, some interesting and some too artificial to be interesting.
\medskip

Let $M$ be a riemannian symmetric, let $\sigma$ 
be an involutive isometry of $M$, let $S$ be a totally geodesic submanifold
of $M$, and suppose that (i) $S$ is a topological component of the fixed 
point set of $\sigma$, (ii) $\dim_\R S = \frac{1}{2}\dim_\R M$, and 
(iii) $S$ has quaternionic 
structure for which its holonomy has quaternion scalar part.  Then we will
say that $S$ is a {\bf quaternionic form} of $M$.
\medskip

Suppose that $M = G/K$ with base point $x_0 = 1K$ and $S = L(x_0) = L/V$
where $L$ is the identity component of the group generated by transvections
of $S$.  Following Proposition \ref{quaternion-criterion}, $S = L/V$ is
one of the spaces listed in Table \ref{quaternion-symm}.  That 
gives us interesting examples 
\medskip

$SU(r+2)/S(U(r) \times U(2)) = U(r+2)/[U(r) \times U(2)]$ in 
$Sp(r+2)/[Sp(r) \times Sp(2)]$;
\smallskip

$SO(r+4)/[SO(r) \times SO(4)]$ in $U(r+4)/[U(r) \times U(4)]
= SU(r+4)/S(U(r) \times SU(4))$;
\smallskip

$Sp(r+1)/[Sp(r) \times Sp(1)]$ in $U(2r+2)/[U(2r) \times U(2)] 
= SU(2r+2)/S(U(2r) \times U(2))$;
\smallskip

$SU(r+2)/S(U(r) \times U(2)) = U(r+2)/[U(r) \times U(2)]$ in
$SO(2r+4)/[SO(2r) \times SO(4)]$;

$P^2(\H) = Sp(3)/[Sp(2) \times Sp(1)]$ in $P^2({\mathbb O}) = F_4/Spin(9)$
\phantom{XXXi} (computing with LiE);
\smallskip

$E_7/[Spin(12)\cdot Sp(1)]$ in $E_8/SO(16)$  (using Proposition
\ref{equal-rank-exist} with $L = E_7A_1$\,, as in \S \ref{equal-rank-exceptional}).
\medskip

\noindent It also gives us some other examples 
\medskip

$S$ in $S \times S$ as a factor or as the diagonal;
\smallskip

$SU(r+2)/S(U(r) \times U(2))$ in $SU(r+4)/S(U(r) \times U(4))$
or $SU(2r+2)/S(U(2r) \times U(2))$;
\smallskip

$SO(r+4)/[SO(r) \times SO(4)]$ in $SO(r+8)/[SO(r) \times SO(8)]$
or $SO(r+8)/[SO(r) \times SO(8)]$;
\smallskip

$Sp(r+1)/[Sp(r) \times Sp(1)]$ in $Sp(r+2)/[Sp(r) \times Sp(2)]$
or $Sp(2r+1)/[Sp(2r) \times Sp(1)]$.
\medskip

\noindent Those other examples somehow seem too formal to be interesting.
Of course with any of these compact examples $S \subset M$, we also have the
noncompact duals $S' \subset M'$.
\medskip

These examples indicate that a reasonable theory for quaternionic 
forms $S$ of symmetric spaces $M$ will require some additional structure 
on the normal bundle of $S$ in $M$.

\enddocument
\end